\documentclass[11pt]{article}
\usepackage{amsmath,amsfonts,amssymb,latexsym,amsbsy, bbm, theorem,enumerate,color}
\usepackage{graphicx, graphics}
\usepackage{float,color,fancybox,shapepar,setspace,hyperref}
\usepackage[affil-it]{authblk}
\usepackage{caption}
\usepackage{tikz}
\usetikzlibrary{arrows}
\tikzset{arr/.style={line width=.5mm, {-Latex[left]}, #1}}
\usepackage[normalem]{ulem}
\usetikzlibrary{arrows,automata}
\usepackage[latin1]{inputenc}
\usepackage{verbatim}

\usepackage{subfigure}
\usepackage{psfrag}
\usepackage{epstopdf}
\newtheorem{thm}{Theorem}[section]
\newtheorem{conj}[thm]{Conjecture}
\newtheorem{lem}[thm]{Lemma}%[section]
\theorembodyfont{\rmfamily}
\def\pf{\medskip\noindent {\bf Proof.}~~}

\def\dfn#1{{\sl #1}}

\newcommand{\less}{\setminus}
\def\qed{ \hfill $\blacksquare$}

\begin{document}

\title{Antimagic orientations of graphs with large maximum degree}

\author{Donglei Yang$^{1}$,  Joshua Carlson$^2$, Andrew Owens$^3$,  K. E. Perry$^4$\\ Inne Singgih$^5$,  Zi-Xia Song$^{6,}$\thanks{Partially supported by the National   Science  Foundation under Grant No. DMS-1854903.}\,,  Fangfang Zhang$^{7,}$\thanks{Corresponding author. This work was done in part while   Fangfang Zhang visited the University of Central Florida as a visiting student.   The visit was 
supported by the Chinese Scholarship Council.    E-mail addresses: dlyang120@163.com (D. Yang);     jc31@williams.edu (J. Carlson);   adowens@widener.edu (A.  Owens);  kperry@soka.edu (K. E. Perry);  Zixia.Song@ucf.edu (Z-X. Song);  Fangfangzh@smail.nju.edu.cn (F. Zhang);  zhangx42@myumanitoba.ca (X. Zhang).}\,,  Xiaohong Zhang$^8$}

  \affil{
  { \small {$^1$Department  of Mathematics, Shandong University, Jinan  250100, China}}\\
   { \small {$^2$Department of Mathematics and Statistics, Williams College, Williamstown, MA 01267, USA}}\\
    { \small {$^3$Department of Mathematics, Widener University, Chester, PA 19013, USA}}\\
     {\small {$^4$Mathematics, Soka University of America, Aliso Viejo, CA 92656}}\\
  { \small {$^5$Department of Mathematics, University of South Carolina, Columbia, SC 29208, USA}}\\
  { \small {$^6$Department  of Mathematics, University of Central Florida, Orlando, FL 32816, USA}} \\
   { \small {$^7$Department of Mathematics, Nanjing University, Nanjing  210093, China}}\\
  { \small {$^8$Department of Mathematics, University of Manitoba, Winnipeg, MB, Canada R3T 2N2}}\\
  }

 \date{ }

\maketitle
\begin{abstract}
Given a digraph $D$ with $m $ arcs, a bijection $\tau: A(D)\rightarrow \{1, 2, \ldots, m\}$ is     an \dfn{antimagic labeling} of $D$ if no two vertices in $D$ have the same  vertex-sum,  where the  vertex-sum  of a vertex $u $  in $D$ under  $\tau$ is the sum of labels of all arcs entering $u$ minus the sum of labels of all arcs leaving $u$.  We say $(D, \tau)$ is  an
  \dfn{antimagic orientation}   of a graph $G$  if    $D$  is an orientation of $G$ and $\tau$ is       an antimagic labeling of $D$.
Motivated by the conjecture of Hartsfield and Ringel from  1990 on antimagic labelings of graphs, Hefetz, M\"{u}tze, and Schwartz in 2010 initiated the study of antimagic  orientations of graphs, and  conjectured that every connected graph  admits an antimagic orientation.   This conjecture seems  hard, and few related results are known. However, it has been verified
to be true for  regular graphs and  biregular bipartite graphs.
 In this paper,  we prove that every connected graph $G$ on $n\ge9$ vertices with maximum degree at least $n-5$  admits an antimagic orientation.

\end{abstract}

\baselineskip 18pt
\section{Introduction}

In this paper,  all  graphs   are finite and  simple, and all multigraphs are finite and loopless.     For a graph $G$, we use $V(G)$, $E(G)$, $|G|$,  $e(G)$, $\Delta(G)$ and $\delta(G)$ to denote the vertex set, edge set, number of vertices, number of   edges, maximum degree, and minimum degree of $G$, respectively.  We use the convention   ``$S:=$'' to mean that $S$ is defined to be the right-hand side of the relation.   For a positive integer $k$, we define $[k]:=\{1,2, \ldots, k\}$.  We use $P_n$,  $C_n$, $K_n$, $S_{n-1}$  and $W_{n-1}$ to denote the path,    cycle, complete graph, star and   wheel  on $n$ vertices, respectively. 

\medskip

  An \dfn{antimagic labeling} of a graph $G$ with $m$ edges is a bijection $\tau:  E(G) \rightarrow  [m]$ such that for any distinct vertices $u$ and $v$, the sum of labels on edges incident to $u$ differs from that for edges incident to $v$.   A graph   is \dfn{antimagic} if it has an antimagic labeling. Hartsfield and Ringel~\cite{NG} introduced antimagic labelings in 1990 and conjectured that every connected graph  other than $K_2$  is  antimagic.  The most recent  progress on this problem is a  result of   Eccles~\cite{E}, which states that   there exists an absolute constant $c_0$ such  that if   $G$  is a graph with average
degree at least $c_0$, and $G$ contains no isolated edge and at most one isolated vertex, then $G$ is
antimagic.  This    improves  a  result of Alon, Kaplan, Lev, Roditty, and Yuster~\cite{Alon}, which states that there exists an absolute   constant  $c$ such that  every   graph on $n$ vertices with minimum degree at least $c\log n$  is antimagic.  Hartsfield and Ringel's Conjecture has also been verified to be true for $d$-regular graphs with $d\ge2$ (see~\cite{DC, DYZ,  KAM, CLPZ}),  and   graphs $G$  with $\Delta(G)\geq|G|-3\geq6$ by   Yilma \cite{Yilma}. For more information on antimagic labelings of graphs and related labeling problems, see the  recent  informative survey~\cite{JAG}.\medskip

Motivated by Hartsfield and Ringel's Conjecture, Hefetz, M\"{u}tze, and Schwartz~\cite{DTJ} introduced antimagic labeling  of digraphs. Let $D$ be a  digraph  and $A(D)$ and $V(D)$ be the set of arcs and vertices of $D$, respectively.  Let $\{a_1, \ldots, a_p\}$ be a set of $p\ge m:=|A(D)|$ positive integers. For an injective mapping    $\tau: A(D)\rightarrow \{a_1,   \ldots, a_p\}$ and for all $u\in V(D)$, we define $s_{(D,\tau)}(u)$  to be  the sum of labels of all arcs entering $u$ minus the sum of labels of all arcs leaving $u$ under $\tau$ when $u$ is not an isolated vertex in $D$, and let $s_{(D,\tau)}(u):=0$ when $u$ is an isolated vertex in $D$.   A bijection $\tau: A(D)\rightarrow [m]$ is     an \dfn{antimagic labeling} of $D$ if  $s_{(D,\tau)}(u)\neq s_{(D,\tau)}(v)$ for all   distinct vertices $u$ and $v$ in $D$.     A digraph $D$ is \dfn{antimagic} if it has an antimagic labeling.  We say $(D, \tau)$ is  an  \dfn{antimagic orientation}   of a graph $G$  if    $D$  is an orientation of $G$ and $\tau$ is       an antimagic labeling of $D$. 
 Hefetz, M\"{u}tze, and Schwartz~\cite{DTJ}  proved that   every orientation of       $S_n$ with $n\neq2$,  $K_n$ with $n\ne3$, and   $W_n$ with $n\ge3$ is antimagic, they further 
 asked whether it is true that    every orientation of any connected  graph, other than $K_3$ and $P_3$, is antimagic.   The same authors   proved an analogous result of Alon, Kaplan, Lev, Roditty, and Yuster~\cite{Alon}, which states that there exists an absolute constant $c$ such that every orientation of any  graph on $n$ vertices with minimum degree at least $c\log n$ is antimagic.  As pointed out in \cite{DTJ}, ``Proving that every orientation of such a graph is antimagic, however, seems rather difficult.''  As a relaxation of this problem, Hefetz, M\"{u}tze, and Schwartz~\cite{DTJ}  proposed the following conjecture.
 
 \begin{conj}[\cite{DTJ}]\label{conj}
Every connected graph admits an antimagic orientation.
\end{conj}
Very recently, Conjecture~\ref{conj} has been verified to be true for  regular graphs   (see \cite{DTJ,   LSWY, Y}),  and biregular bipartite graphs with minimum degree at least two by  Shan and Yu~\cite{SY}. 
Motivated by the work of Yilma~\cite{Yilma},   we establish more evidence for  Conjecture \ref{conj} in this paper by studying antimagic orientations of   graphs with large maximum degrees. We prove the following main result.

\begin{thm}\label{main} Let $G$ be a connected graph. 
\begin{enumerate}[\rm(i)] 
\item If  $\Delta(G)\geq|G|-3$, then $G$ has an antimagic orientation.
\item If  $\Delta(G)=|G|-t\ge4$ for   each  $t\in\{4,5\}$,  then   $G$ has an antimagic orientation.
\end{enumerate} 
\end{thm}
 \medskip
 
We prove Theorem~\ref{main} in Section~\ref{Main}.  The main idea of the proofs  of  Theorem~\ref{main}  and the preliminary results (stated and proved in Section~\ref{preliminary})  is to use Eulerian orientations. This strategy  was  first used   in      \cite{DTJ} and then in \cite{LSWY, Y}. Our method here is more involved to   obtain        antimagic orientations of such graphs.\medskip

We conclude this section by introducing more notation.  
  Given a graph $G$,    sets $S\subseteq V(G)$ and $F\subseteq E(G)$,  we use        $G\less S$ to denote the subgraph    obtained from $G$ by deleting all vertices in $S$,   $G\less F$ the subgraph    obtained from $G$ by deleting all edges in $F$, and $G[S] $    the  subgraph    obtained from $G$ by deleting all vertices in $V(G)\less S$.  We simply write $G\less v$ when $S=\{v\}$.         For  two disjoint sets $A, B\subseteq V(G)$,    $A$ is \dfn{complete} to $B$ in $G$  if each vertex in $A$ is adjacent to all vertices in  $B$.   
We simply say $a$ is \dfn{complete to}  $B$ when $A=\{a\}$.
For convenience, we use   $A \less B$ to denote  $A -B$; and $A \less b$ to denote  $A -\{b\}$ when $B=\{b\}$.  For each vertex $v\in V(G)$, let $N(v)$ denote the set of all vertices adjacent to $v$ in $G$ and let $N[v]:=N(v)\cup\{v\}$.    A closed walk in a multigraph is an \dfn{Euler tour} if
it traverses every edge of the graph exactly once. The following is a result of Euler.
\begin{thm}[Euler 1736]\label{Euler}
A connected multigraph admits an  Euler tour if and only if every vertex has even degree.
\end{thm}
\medskip

\section{Preliminaries}\label{preliminary}
   We begin this section with Lemma~\ref{euler} which will be applied  in the proofs of the remaining results in this paper.

\begin{lem}\label{euler}
Let $G$ be a graph with $m\ge1 $ edges and  let  
 $a_1, \ldots,a_m$ be $m$   positive  integers with $  a_1< \cdots<a_m$.   Then  there exist  an orientation $D$ of $G$ and a bijection $\tau: A(D)\to \{a_1,\ldots,a_m\}$ such that for all     $v\in V(G)$,
\[s_{(D, \tau)}(v)\ge -a_m.\]
\end{lem}
\pf  Let $G$ and $a_1, \ldots, a_m$ be given as  in the statement.   We may assume that          $G$ is connected.  Let $A$ be the set (possibly empty) of all vertices $v\in V(G)$ with $d(v)$ odd. Then $|A|=2t$ for some integer $t\ge0$. Let $G^*:=G$ when $t=0$. When $t\ge1$, we may assume that $A:=\{x_1, x_2, \ldots, x_{2t}\}$.  Let $G^*$ be   obtained from $G$ by adding $t$ new edges $x_ix_{i+t} $ for all $i\in[t]$.  Then $e(G^*)=m+t$. By Theorem~\ref{Euler},  $G^*$ contains an Euler tour  with vertices and edges $v_1, e_1, v_2, e_2, \ldots v_{m+t}, e_{m+t}, v_1$ in order, where $v_1, \ldots, v_{m+t}$ are not necessarily distinct, and edges $e_1=v_1v_2, e_2=v_2v_3, \ldots, e_{m+t}=v_{m+t}v_1$ are pairwise distinct. We may further assume that  $e_1\in E(G)$.  Let $1=i_1<i_2<\cdots< i_m\le m+t$ be such that $ e_{i_1},\ldots, e_{i_m}$ are all the edges of $G$.   Let $D$ be the orientation of $G$ obtained by orienting  each edge $e_{i_j}$ from $v_{i_j}$ to $v_{{i_j}+1}$ for all $j\in[m]$.  Let $\tau: E(G)\rightarrow \{a_1, \ldots, a_m\}$ be the bijection such that $\tau(e_{i_j})=a_{m+1-j}$ for all $j\in[m]$.  It can be   checked that for all $v\in V(G)$,   $s_{(D, \tau)}(v)\ge -a_m$.   \medskip

This completes the proof of Lemma~\ref{euler}.
\qed\bigskip

\noindent {\bf Remark}.  For the sake of simplicity and clarity of presentation, we  shall  apply  Lemma~\ref{euler} to graphs $H$ with no edges. Under those circumstances, we  shall   let   $D$ with    $V(D)=V(H)$ and   $A(D)=\emptyset$ be the orientation of $H$,   and $\tau: A(D)\rightarrow   \emptyset$  with $s_{(D, \tau)}(v)=0$ for all $v\in V(D)$  be the bijection. \medskip 

 \begin{lem}\label{independent}
Let $G$ be a connected graph and let $x \in V(G)$ be a vertex of degree     $\Delta(G)\le |G|-2$.   If $V(G)\less N[x]$ is an independent set, then $G$ admits an antimagic orientation.
\end{lem}
\pf Let $G$ and $x$ be given as in the statement. Then $t:=|V(G)\less N[x] |\ge1 $ because $d(x)\le |G|-2$.
Let $V(G)\less N[x]:=\{y_1, \ldots, y_{t}\}$. Let $d:=d(x)$ and $d_i:=d(y_i)$ for all $i\in[t]$.   Since $G$ is connected, we may assume that  $1\le d_1\le\cdots\le d_t $.   Let $m_1:=  d_1+\cdots+d_t$ and let $e_1,  \ldots, e_{m_1}\in E(G)$ be  all the  edges between $\{y_1, \ldots, y_{t}\}$ and $N(x)$ in $G$ such that $y_1$ is incident with $e_1, \ldots, e_{d_1}$, and for all $i\in\{2, \ldots, t\}$, $y_i$ is incident with  $e_{d_1+\cdots +d_{i-1}+1}, \ldots, e_{d_1+\cdots +d_{i}}$. By Lemma~\ref{euler} applied to $H:=G\less\{x, y_1, \ldots, y_t\}$, there exist a bijection     $\tau': E(H)\rightarrow  \{m_1+1, \ldots, m-d \}$  and an orientation $D'$  of  $H$ such that  $s_{(D', \tau')}(v)\geq -(m-d )$ for all $v\in N(x)$.
Let $D$ be the orientation of $G$ obtained from $D'$ by orienting all the edges  between $\{x, y_1, \ldots, y_t\}$ and $N(x)$ towards $N(x)$. Let $\tau'': A(D\less x)\rightarrow   [m-d ]$ be the bijection obtained from $\tau'$ by letting $\tau''(e_i)=i$ for all $i\in[m_1]$. We may further assume that $N(x)  :=\{x_1,  \ldots, x_{d }\}$ with  $s_{(D\less x, \tau'')}(x_1)\le \cdots\le  s_{(D\less x, \tau'')}(x_d)$.   Finally, let  $\tau: A(D)\rightarrow   [m]$ be the bijection obtained from $\tau''$ by letting $\tau(xx_i)=m-d+i$ for all $i\in[d]$.  Then for all   $i\in [d]$,  \[s_{(D, \tau)}(x_i)=s_{(D\less x, \tau'')}(x_i)+(m-d+i)\ge  -(m-d)+ (m-d +i)>0. \]  By the choice of $(D, \tau)$ and the fact that $d(x)=\Delta(G)$, it follows that
\[s_{(D, \tau)}(x)<s_{(D, \tau)}(y_t) <\cdots<s_{(D, \tau)}(y_1)<0<s_{(D, \tau)}(x_1)  <\cdots<s_{(D, \tau)}(x_d).\]
Hence $(D, \tau)$ is an  antimagic orientation of $G$.
This completes the proof of Lemma~\ref{independent}.
\qed\bigskip

We end this section by proving that Conjecture~\ref{conj} is true for graphs $G$ with    $|G|\ge9$  and     a dominating set of size two.  Theorem~\ref{dominating} implies that every complete multipartite graph admits an antimagic orientation. \medskip
 
\begin{thm}\label{dominating}
Let $G$ be a graph and let  $x,y\in V(G)$ be distinct  such that   $N[x]\cup N[y] =V(G)$ and $d(x)\ge d(y)$.  If     $d(x)\ge4$  or  $ N(x)\cap N[y] \ne\emptyset $, then $G$ admits an antimagic orientation.
\end{thm}

\pf  Let $G$ and $x,y$ be given as in the statement.  Then $d(x)\ge d(y)$. Let $n:=|G|$, $m:=e(G)$ and $A:= V(G)\less\{x,y\}$.   Since $N[x]\cup N[y] =V(G)$,  we see that every vertex in $A$ is adjacent to either $x$ or $y$ in $G$.     \medskip

Assume first that $   N(x)\cap N[y] =\emptyset $.      Then   $d(x)\ge4$, $xy\notin E(G)$, $N(x)\cap N(y)=\emptyset$ and $N(x)\cup N(y)=A$.            Let $e_v$ be the unique edge between $v $ and $\{x,y\}$ in $G$ for all $v\in A$. By  Lemma~\ref{euler} applied to $H:=G\less\{x,y\}$, there exist an orientation $D'$  of  $H$ and a bijection     $\tau': A(D')\rightarrow  [m-n+2]$  such that   $s_{(D',\tau')}(v)\geq -(m-n+2)$ for all   $ v\in A$.
   We may further assume that $A   =\{v_1,v_2,\ldots, v_{n-2}\}$ with  $s_{(D',\tau')}(v_1)\le \cdots\le  s_{(D', \tau')}(v_{n-2})$.    Let $D$ be the orientation of $G$ obtained from $D'$ by orienting each edge $e_{v_i}$ towards $v_i $ for all $i\in[n-2]$. Let $\tau: A(D)\rightarrow [m]$ be the  bijection obtained from   $\tau'$ by letting $\tau(e_{v_i})=m-n+2+i$ for all $i\in [n-2]$.
Then  by the choice of $(D, \tau)$,
$0<s_{(D, \tau)}(v_1)<\cdots<s_{(D, \tau)}(v_{n-2}) $,   because  for all  $ i\in [n-2]$,
\[s_{(D, \tau)}(v_i)=s_{(D',\tau')}(v_i)+(m-n+2+i)\ge  -(m-n+2)+ (m-n+2+i)>0. \] Note that $s_{(D, \tau)}(x)<0$ and $s_{(D, \tau)}(y)\le0$.  Then $(D, \tau)$ is an  antimagic orientation of $G$ if  $s_{(D, \tau)}(x)  \ne s_{(D, \tau)}(y) $. We may assume that    $s_{(D, \tau)}(x)  = s_{(D, \tau)}(y) $.   Since  $d(x)\ge4$, we see that    $s_{(D, \tau)}(x)<-4(m-n+3)$. By the choice of $\tau$, $\tau(e_{v_1})=m-n+3$.  Let $x', y'$ be a permutation of $x, y$ such that  $e_{v_1}=v_1x'$. Let $ D^* $ be obtained from $ D $ by   reorienting the edge $v_1x'$ away from the vertex $v_1$. Then   
$s_{(D^*, \tau)}(x')=s_{(D, \tau)}(x) +2(m-n+3) <-2(m-n+3)$, $s_{(D^*, \tau)}(y')=s_{(D, \tau)}(x)<s_{(D^*, \tau)}(x')$, $s_{(D^*, \tau)}(v_1)=s_{(D, \tau) }(v_1)-2(m-n+3) \ge  -(m-n+2) -(m-n+3)=-2(m-n+3)+1$,  and $s_{(D^*, \tau)}(v)= s_{(D, \tau)}(v)$ for all $v\in V(G)\less \{v_1, x'\}$. It follows that
\[s_{(D^*, \tau)}(y')<s_{(D^*, \tau)}(x') <s_{(D^*, \tau)}(v_1)  <\cdots<s_{(D^*, \tau)}(v_{n-2}).\]
Hence $(D^*, \tau)$ is an  antimagic orientation of $G$. \medskip

Assume next that   $t:= |N(x)\cap N[y]| \ge1$.   Then $d(y)\ge1$. Let $N(x)\cap N[y]:=\{x_1, \ldots, x_t\}$ such that $x_1=y$ if $xy\in E(G)$.  Then $x$ is complete to $\{x_1, \ldots, x_t\}$  and $y$ is complete to $\{x_1, \ldots, x_t\}\less \{y\}$ in $G$.   The statement holds trivially when $n\le 3$. We may assume that $n\ge4$. For all $v\in A$, let $e_v$ be the unique edge between $v $ and $\{x,y\}$ in $G\less\{xx_1, \ldots, xx_t\}$.  By  Lemma~\ref{euler} applied to $ H:=G\less\{x,y\}$, there exist a bijection     $\tau': E(H)\rightarrow \{t+1, \ldots, m-n+2 \}$  and an orientation $D'$  of  $ H$  such that $s_{(D', \tau')}(v)\geq -(m-n+2)$ for all $v\in A$.   Let $D''$ be the orientation of $G \less \{e_v: v\in A\}$ obtained from $D'$ by orienting the edge $xx_1$ away from the vertex $y$ when $x_1=y$ and   every edge  $xx_i$  away from the vertex $x$ for all $i\in[t]$ with $x_i\ne y$.   Let $\tau'': E(G)\less  \{e_v: v\in A\}\rightarrow  [m-n+2]$  be the bijection obtained from $\tau'$ by letting    $\tau''(xx_j)=j$ for all $j\in [t]$.
We may further assume that $A   =\{v_1,v_2,\ldots, v_{n-2}\}$ with  $s_{(D'',\tau'')}(v_1)\le \cdots\le  s_{(D'', \tau'')}(v_{n-2})$.
Finally, let $D$ be the orientation of $G $ obtained from $D''$ by orienting each edge $e_v$ towards the vertex $v$ for all $v\in A$, and  let $\tau: A(D)\rightarrow [m]$ be the  bijection obtained from   $\tau''$ by letting $\tau(e_{v_i})=m-n+2+i$ for all $i\in [n-2]$.
Then  by the choice of $(D, \tau)$,
$0<s_{(D, \tau)}(v_1)<\cdots<s_{(D, \tau)}(v_{n-2}) $,   because $s_{(D, \tau)}(v_i)\ge s_{(D'', \tau'')}(v_i)+(m-n+2+i)\ge  -(m-n+2)+ (m-n+2+i)>0$ for all $i\in [n-2]$.  Note that $s_{(D, \tau)}(x)<0$ and $s_{(D, \tau)}(y)<0$.  Then $(D, \tau)$ is an  antimagic orientation of $G$ if $s_{(D, \tau)}(x)  \ne s_{(D, \tau)}(y) $. So we may assume that    $s_{(D, \tau)}(x)  = s_{(D, \tau)}(y) $.     Then $d(y)>1$ because $n\ge4$.  Assume that      $x_i$ is complete to $\{x, y\}$ in $G$ for some $i\in[t]$.   Let $\tau^*: A(D)\rightarrow [m]$ be the  bijection obtained from   $\tau$ by letting $\tau^*(xx_i)=\tau(yx_i)$, $\tau^*(yx_i)=\tau(xx_i)$ and $\tau^*(e)=\tau(e)$ for all $e\in E(G)$ with $e\ne xx_i, yx_i$. Then
 $(D, \tau^*)$ is an  antimagic orientation of $G$.  Next,  assume that    $t=1$ and $xy\in E(G)$.   Let $D^*$ be the orientation of $G$ obtained from $D$ by reorienting the edge $xy$ towards the vertex $y$. Then $(D^*, \tau)$ is an  antimagic orientation of $G$.\medskip

 This completes the proof of Theorem \ref{dominating}.
 \qed\medskip

\section{Proof of Theorem \ref{main}}\label{Main}

We are now ready to prove Theorem \ref{main}.  
Let  $G$ be a connected graph.    Let  $x\in V(G)$ be  such that  $d(x)=\Delta(G)$. Let $Y:=V(G)\less N[x]$. To prove Theorem \ref{main}(i), assume that $d(x)\ge |G|-3$. 
By Lemma~\ref{independent}, we may assume that $d(x)=|G|-1$ or $d(x)=|G|-3$ and $G[Y]=K_2$.   Let $y\in N(x)$ be an arbitrary vertex when $d(x)=|G|-1$,  and let $y\in Y$   such that $N(y)\cap N(x)\ne\emptyset$  when $\Delta(G)=|G|-3$ (this is possible because $G$ is connected).  Then $N[x]\cup N[y]=V(G)$ and $  N(x)\cap N[y]\ne\emptyset$.   By Theorem~\ref{dominating}, $G$ admits an antimagic orientation. \medskip

  To prove Theorem \ref{main}(ii), assume that $d(x)= |G|-t\ge4$, where $t\in\{4,5\}$.   Then $|G|\ge t+4$ and $d(x)\ge4$. Let   $ Y:=\{y_1,\ldots,y_{t-1}\}$.    By Lemma~\ref{independent} and Theorem~\ref{dominating}, we may assume that   $1\le \Delta(G[Y])\le t-3$.  Let $Y_1, \ldots, Y_r$ be all components of $G[Y]$  with $|Y_1|\ge  \cdots\ge|Y_r|$,
where $1\le r\le t-2$.  Choose  $z_1, \ldots, z_r\in N(x)$ (not necessarily distinct) such that for each $i\in[r]$, $z_i$ is adjacent to a vertex in $Y_i$. This is possible because $G$ is connected. For each $i\in[r]$, let $e_i$ be an arbitrary edge between  $z_i$ and $Y_i$ in $G$. Let $H_0$ be the subgraph of $G$ with $V(H_0)=Y\cup\{z_1, \ldots, z_r\}$ and $E(H_0)=\{e_i: i\in [r]\}\cup E(G[Y])$.   Let $H_1$ be the bipartite subgraph of $G$ with bipartition $\{Y, N(x)\}$ and $E(H_1) =E_G(Y, N(x))\less \{e_i: i\in [r]\}$, where  $E_G(Y, N(x))$ denotes    the set of   all edges between $Y$ and $N(x)$ in $G$. \medskip

\begin{figure}[htbp]
\centering
\includegraphics[scale=0.8]{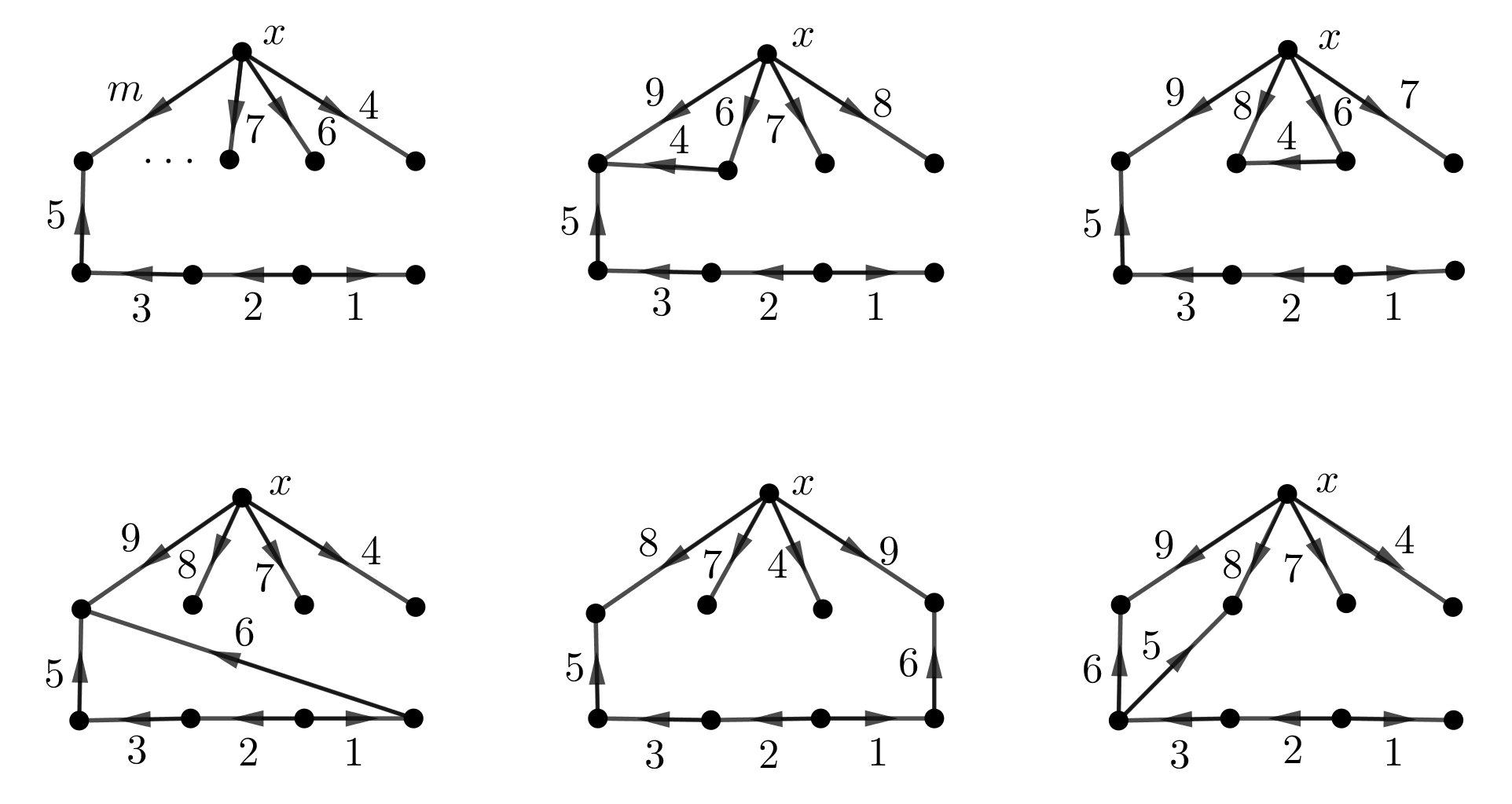}
\caption{Antimagic orientations of all possible $G$ when $H_0=P_5$ and $m\in\{8, 9\}$.}
\label{m89}
\end{figure}
\begin{figure}[htbp]
\centering
\includegraphics[scale=0.8]{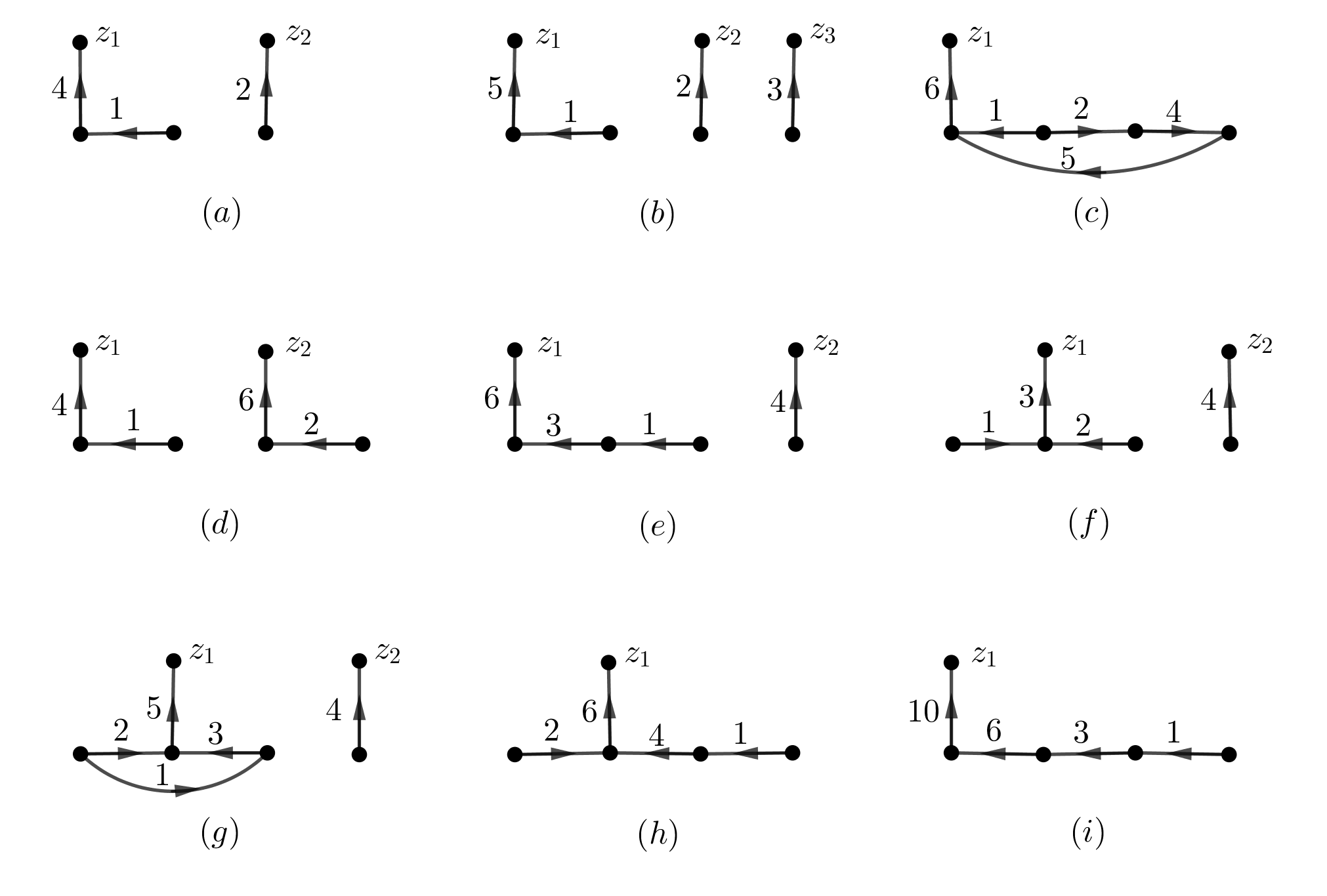}
\caption{$(D_0,\tau_0)$  of  all possible $H_0$ according to  $e(G[Y])$.}
\label{D0}
\end{figure}

Let  $m:=e(G)$. Then $m\ge |G|-1\ge t+3$.  Note that if $G[Y]=P_4$, then $Y_1=G[Y]$, $t=5$ and $m\ge t+3 = 8$. In this case,  we choose $z_1$ so that  $ z_1$ is adjacent to a vertex of degree two in $G[Y]$ if possible.  Hence, the statement holds when $H_0=P_5$ and $m\le 9$, as depicted in Figure~\ref{m89}.  We may assume that $m\ge10$ when $H_0=P_5$.  We next find an orientation $D_0$ of $H_0$  and an injective mapping  $\tau_0: A(D_0)\rightarrow \{1,\ldots,6\}$ when $H_0\ne P_5$,   and $\tau_0: A(D_0)\rightarrow \{1,\ldots,6,10\}$ when  $H_0= P_5$  and $m\ge10$  such that $ -4\le  s_{(D_0, \tau_0)}(y_i)\le 0$ for all $i\in [t-1]$,    and  $0<|s_{(D_0, \tau_0)}(y_i)- s_{(D_0, \tau_0)}(y_j)|\le 4$ for $1\le i<j\le t-1$.  Then $s_{(D_0, \tau_0)}(y_1), \ldots, s_{(D_0, \tau_0)}(y_{t-1})$  are pairwise distinct.  Such a $(D_0, \tau_0)$  is depicted in  Figure~\ref{D0}(a,b) when $e(G[Y])=1$, Figure~\ref{D0}(d,e,f) when $e(G[Y])=2$, Figure~\ref{D0}(g,h,i) when $e(G[Y])=3$, and Figure~\ref{D0}(c) when $e(G[Y])=4$, respectively.  \medskip

For $i, j\in [t-1]$ with $i\ne j$, we define  $y_i\prec y_j$ if either $d_{H_1}(y_i)<d_{H_1}(y_j)$, or $d_{H_1}(y_i)=d_{H_1}(y_j)$ and $s_{(D_0, \tau_0)}(y_i)>s_{(D_0, \tau_0)}(y_j)$. This is possible because  $s_{(D_0, \tau_0)}(y_1), \ldots, s_{(D_0, \tau_0)}(y_{t-1})$  are pairwise distinct.   We may assume that  $y_1\prec y_2\prec \cdots\prec y_{t-1}$.  Let $N_0:=\tau_0(A(D_0))$, $N_1$ be the set (possibly empty) of $d_{H_1}(y_1)$ many smallest numbers  in   $\{5,\ldots,m\}\less N_0$, and for $i\in\{2, \ldots, t-1\}$, let $N_i$ be the set (possibly empty) of  $d_{H_1}(y_i)$ many smallest numbers  in   $\{5,\ldots,m\}\less (N_0\cup N_1\cup\cdots\cup N_{i-1})$.   Let $D_1$ be an orientation of $H_0\cup H_1$ obtained from $D_0$ by orienting all edges in $H_1$ away from $Y$, and $\tau_1:E(H_0)\cup E(H_1)\rightarrow N_0\cup N_1\cup\cdots\cup N_{t-1}$ be the bijection obtained from $\tau_0$ by first letting $\tau_1(e)=\tau_0(e)$ for all $e\in E(H_0)$ and then assigning all the numbers in $N_i $  to the edges incident with $y_i$ in $H_1$ for all $i\in[t-1]$.
  By the choice of $(D_1, \tau_1)$ and the fact that  $y_1\prec y_2\prec \cdots\prec y_{t-1}$,  for $i, j\in[t-1]$ with $i< j$,  $ s_{(D_1, \tau_1)}(y_i)- s_{(D_1, \tau_1)}(y_j)  >0$ because  $|s_{(D_0, \tau_0)}(y_i)-s_{(D_0, \tau_0)}(y_j)|\le 4$ if $d_{H_1}(y_i)<d_{H_1}(y_j)$, and  $ s_{(D_0, \tau_0)}(y_i)- s_{(D_0, \tau_0)}(y_j)  >0$  if $d_{H_1}(y_i)=d_{H_1}(y_j)$.  It follows that
  \[0\ge s_{(D_1, \tau_1)}(y_1)>  s_{(D_1, \tau_1)}(y_2)>\cdots>s_{(D_1, \tau_1)}(y_{t-1}).\]

Next, let $H_2:= G[N(x)]$,  $m_1:=e(H_2)$ and let $N$ be the set (possibly empty) of $m_1$ smallest numbers in $[m]\less (N_0\cup N_1\cup\cdots\cup N_{t-1})$.  Let $p:=0$ when $N=\emptyset$,  and $p$ be the largest number in $N$ when $N\ne \emptyset$. By Lemma \ref{euler} applied to $H_2$, there exist an orientation $D_2$ of     $H_2$ and a bijection $\tau_2: A(D_2)\rightarrow  N$ such that $s_{(D_2, \tau_2)}(v)\ge -p $ for all $v\in N(x)$. Let   $n:=|G|$ and $N(x):= \{v_1,v_2,\ldots, v_{n-t}\}$ such that   \[s_{(D_1,\tau_1)}(v_1)+s_{(D_2,\tau_2)}(v_1)\le s_{(D_1,\tau_1)}(v_2)+s_{(D_2,\tau_2)}(v_2)\le \cdots\le  s_{(D_1,\tau_1)}(v_{n-t})+s_{(D_2,\tau_2)}(v_{n-t}).\]
Let $ a_1, \ldots, a_{n-t } $  be the remaining   $n-t $ numbers in $[m]\less (N_0\cup N_1\cup\cdots\cup N_{t-1}\cup N)$.   We may further assume that   $a_1<a_2<\ldots<a_{n-t }$. Then $p<a_1$. Let $D$ be the orientation  of $G$ obtained from $D_1$ and $D_2$ by orienting all the edges incident with $x$ away from  $x$, and let $\tau: A(D)\rightarrow [m]$ be the bijection obtained from $\tau_1$ and $\tau_2$ by letting $\tau(xv_i)=a_i$ for all $i\in [n-t]$, $\tau(e)=\tau_1(e)$ for all $e\in E(H_0)\cup E(H_1)$, and $\tau(e)=\tau_2(e)$ for all $e\in E(H_2)$. \medskip

It remains to  show that  $(D, \tau)$ is an antimagic orientation of $G$. By the choice of $(D,\tau)$, $s_{(D, \tau)}(y_i) = s_{(D_1, \tau_1)}(y_i)$ for all $i\in[t-1]$. Furthermore, for all $v_i\in N(x)$,   \[s_{(D, \tau)}(v_i)=s_{(D_1,\tau_1)}(v_i)+s_{(D_2,\tau_2)}(v_i)+\tau(xv_i)\ge -p+a_i>0, \text{ because } p<a_i.\]
 It follows that
\[ s_{(D, \tau)}(y_{t-1}) <s_{(D, \tau)}(y_{t-2}) <\cdots< s_{(D, \tau)}(y_1)\le 0<s_{(D, \tau)}(v_1)< s_{(D, \tau)}(v_2)<\cdots < s_{(D, \tau)}(v_{n-t}).\]
 We next  show that
  $s_{(D, \tau)}(x)<s_{(D, \tau)}(y_{t-1})$.  Note that $m=d(x)+m_1+e(H_1)+e(H_0)$ and $d(x)=\Delta(G)=n-t\ge4$.    
   By the choice of $(D_0,\tau_0)$ and $(D, \tau)$, $a_1\ge2$,  and $a_2=4$ only when $m_1=0$ and $(D_0, \tau_0)$ is depicted in Figure~\ref{D0}(i). Hence,
  \[s_{(D, \tau)}(x) = -\sum_{i=1}^{n-t} a_i\le \begin{cases}-(a_1+1+2)-(n-t-1)a_2\le -5-(n-t-1)a_2 \,\, &\text{ if } H_0\ne P_5, \\
  -(a_1+a_2+1)-(n-t-2)a_3\le -7-(n-t-2)a_3 \,\,& \text{ if } H_0= P_5.
  \end{cases}
  \] Then $s_{(D, \tau)}(x)<-4\le s_{(D, \tau)}(y_{t-1})$ when $N_{t-1}=\emptyset$. We may assume that $N_{t-1}\ne\emptyset$. Then
  \[|N_{t-1}| =d_{H_1}(y_{t-1})\le d_G(y_{t-1})-1\le \Delta(G)-1=n-t-1,\]
   and  $|N_{t-1}| =n-t-1$  only when $d_G(y_{t-1})=\Delta(G)=n -t$ and $d_{H_0}(y_{t-1})=1$.
   Assume first that $H_0= P_5$ and $|N_{t-1}| =n-t-1$. In this case,  $t=5$, $(D_0, \tau_0)$ is  depicted in Figure~\ref{D0}(i),   $d_{H_0}(y_4)=1$ and  $d_G(y_4)=n-5\ge4$. Let $b_1, b_2, \ldots, b_{n-5}$ be the $n-5$ positive integers assigned to the edges incident with $y_4$ under $\tau$.  We may further assume that $b_1<b_2<\cdots<b_{n-5}$.  Then $b_1=1$.  By the choice of $ \tau$,  $b_2\le8$ when $n=9$ and  $b_2\le 4+(n-5)+1=n $ when $n\ge10$,  and $a_{i}-b_{i}\ge n-7 $ for all $i\in\{3, \ldots, n-5\}$. Then
\begin{align*}
s_{(D, \tau)}(y_t)-s_{(D, \tau)}(x)&= (a_1-b_1)+(a_2-b_2)+\sum_{i=3}^{n-5}(a_i-b_i)\\
&\ge    \begin{cases} (2-1)+(4-n)+3(n-7) & \text{ if } n\ge10 \\
  (2-1)+(4-8)+2(n-7)   & \text{ if }  n=9
  \end{cases}\\
&>0.
\end{align*}
  Assume next that $H_0\ne  P_5$ or  $|N_{t-1}| \le n-t-2$.   Let $q^*$ be the largest number in $N_{t-1}$. 
    By the choice of $(D_0,\tau_0)$ and $N_{t-1}$,
  if $H_0\ne P_5$, then $q^*<a_2$;  if $H_0= P_5$, then $q^*<a_3$. Hence
\begin{align*}
s_{(D, \tau)}(y_{t-1})&=s_{(D_0, \tau_0)}(y_{t-1})-\sum_{q\in N_{t-1}} q \\
 &\ge  \begin{cases}-4-(n-t-1)q^*  & \text{ if } H_0\ne P_5  \\
  -4-(n-t-2)q^*   & \text{ if } H_0= P_5 \, \text{ and } \,  |N_{t-1}| \le n-t-2
  \end{cases}\\
  &> \begin{cases}-4  -(n-t-1)a_2  & \text{ if } H_0\ne P_5  \\
  -4 -(n-t-2)a_3  & \text{ if } H_0= P_5 \, \text{ and } \,  |N_{t-1}| \le n-t-2
    \end{cases}\\
  &> s_{(D, \tau)}(x).
 \end{align*}
 In all cases,   we have $s_{(D, \tau)}(x)<  s_{(D, \tau)}(y_{t-1})$.
It follows that
\[s_{(D, \tau)}(x)<s_{(D, \tau)}(y_{t-1}) < \cdots< s_{(D, \tau)}(y_1)\le 0<s_{(D, \tau)}(v_1)< s_{(D, \tau)}(v_2)<\cdots < s_{(D, \tau)}(v_{n-t}).\]
  Hence  $(D, \tau)$ is an antimagic orientation of $G$.
This completes the proof of Theorem~\ref{main}.\qed\bigskip

\noindent{\bf Acknowledgements}\bigskip 

This work was completed in part at the 2019 Graduate Research Workshop in Combinatorics, which was supported in part by NSF grant \#1923238, NSA grant \#H98230-18-1-0017,  a generous award from the Combinatorics Foundation, and Simons Foundation Collaboration Grants \#426971 (to M. Ferrara), \#316262 (to S. Hartke) and \#315347 (to J. Martin).


\bigskip



\begin{thebibliography}{99}

\bibitem{Alon} N. Alon, G. Kaplan, A. Lev, Y. Roditty, and R. Yuster, Dense graphs are antimagic, J. Graph Theory 47 (2004), 297-309.

\vspace{-0.25cm}
 
 \bibitem{KAM} K. B\'{e}rczi, A. Bern\'{a}th, and M. Vizer, Regular graphs are antimagic, Electron. J. Combin. 22 (2015): Paper 3.34.
  \vspace{-0.25cm}
  \bibitem{CLPZ} F. Chang, Y-C. Liang, Z. Pan, and X. Zhu,  Antimagic labeling of regular graphs,  J. Graph Theory 82 (2016), 339--349.
 \vspace{-0.25cm}

 \bibitem{DC} D. Cranston, Regular bipartite graphs are antimagic, J. Graph Theory 60 (2009) 173--182.
  \vspace{-0.25cm}

\bibitem{DYZ} D. W. Cranston, Y-C. Liang, and X. Zhu, Regular graphs of odd degree are antimagic, J. Graph Theory 80 (2015) 28--33.
 \vspace{-0.25cm}
\bibitem{E} T. Eccles,  Graphs of large linear size are antimagic,  J. Graph Theory 81 (2016) 236--261.
 \vspace{-0.25cm}


\bibitem{JAG} J. A. Gallian, A dynamic survey of graph labeling, Electron J. Combin DS6 (2016).
 \vspace{-0.25cm}
\bibitem{NG} N. Hartsfield and G. Ringel, Pearls in Graph Theory, Academic Press, Boston, (1990) 108-109 (revised version, 1994).
\vspace{-0.25cm}
 
 \bibitem{DTJ} D. Hefetz, T. M\"{u}tze, and J. Schwartz, On antimagic directed graphs, J. Graph Theory 64 (2010) 219--232.
 \vspace{-0.25cm}
 
 
\bibitem{LSWY} T. Li, Z-X. Song, G. Wang,   D. Yang, and C-Q. Zhang, Antimagic orientations of even regular graphs, J. Graph Theory 90 (2019) 46--53.
 \vspace{-0.25cm}
  
\bibitem{SY} S. Shan and X. Yu, Antimagic orientation of biregular bipartite graphs, Electron. J. Combin.  24 (2017): Paper 4.31.
\vspace{-0.25cm}

\bibitem{Y} D. Yang, A note on antimagic orientations of even regular graphs, to appear in Discrete Appl.  Math. (https://doi.org/10.1016/j.dam.2019.04.017). 
\vspace{-0.25cm}
\bibitem{Yilma} Z. B. Yilma, Antimagic properties of graphs with large maximum degree, J. Graph Theory 72 (2013) 367--373.

\end{thebibliography}
\end{document}